\documentclass[letterpaper, 10 pt, conference]{ieeeconf}  % Comment this line out if you need a4paper

\IEEEoverridecommandlockouts                              % This command is only needed if 
\usepackage{amsmath} % assumes amsmath package installed
\usepackage{amssymb}  % assumes amsmath package installed
\usepackage{mathtools}
\let\labelindent\relax
\usepackage{enumitem}

\newcommand{\ie}{i.\,e.~}
\newcommand{\eg}{e.\,g.~}
\newcommand{\Span}[1]{\mathrm{span}\{#1\}}
\newcommand{\Rank}[1]{\mathrm{rank}(#1)}
\newcommand{\Dim}[1]{\mathrm{dim}(#1)}
\newcommand{\C}[1]{\mathcal{C}(#1)}

\newcommand{\Lie}{\mathrm{L}}

\newcommand{\Mod}{\mathrm{~mod~}}

\newtheorem{theorem}{Theorem}
\newtheorem{lemma}{Lemma}

\newtheorem{remark}{Remark}

\title{\LARGE \bf A Finite Test for the Linearizability of Two-Input Systems by a Two-Dimensional Endogenous Dynamic Feedback
}

\author{Conrad Gst{\"o}ttner$^{1}$ Bernd Kolar$^{1}$ Markus Sch{\"o}berl$^{1}$% <-this % stops a space
\thanks{*The first author and the second author have been supported by the Austrian Science Fund (FWF) under grant number P 32151 and P 29964.}% <-this % stops a space
\thanks{$^{1}$Institute of Automatic Control and Control Systems Technology,	Johannes Kepler University, Linz, Austria {\tt\small \{conrad.gstoettner,bernd.kolar,\newline \indent\indent markus.schoeberl\}@jku.at}}%
}

\begin{document}

	\maketitle
	\thispagestyle{empty}
	\pagestyle{empty}
	
	%%%%%%%%%%%%%%%%%%%%%%%%%%%%%%%%%%%%%%%%%%%%%%%%%%%%%%%%%%%%%%%%%%%%%%%%%%%%%%%%
	\begin{abstract}
		We propose an algorithmic test to check whether a two-input system is linearizable by an endogenous dynamic feedback with a dimension of at most two. This test furthermore provides a procedure for systematically deriving flat outputs for this class of systems.
	\end{abstract}

	%%%%%%%%%%%%%%%%%%%%%%%%%%%%%%%%%%%%%%%%%%%%%%%%%%%%%%%%%%%%%%%%%%%%%%%%%%%%%%%%
	\vspace{-1.1ex}
	\section{INTRODUCTION}
		\vspace{-0.9ex}
		\noindent
		Flatness was introduced in control theory by Fliess, L\'evine, Martin and Rouchon, see \eg \cite{FliessLevineMartinRouchon:1992,FliessLevineMartinRouchon:1995}. The flatness property allows an elegant systematic solution of feed-forward and feedback problems, see \eg \cite{FliessLevineMartinRouchon:1995}. Roughly speaking, a nonlinear control system
		\vspace{-1.5ex}
		\begin{align*}
			\begin{aligned}
				\dot{x}&=f(x,u)
			\end{aligned}
		\end{align*}
		\vspace{-3.9ex}
		
		\noindent
		with $\Dim{x}=n$ states and $\Dim{u}=m$ inputs is flat, if there exist $m$ differentially independent functions $y^j=\varphi^j(x,u,u_1,\ldots,u_q)$, $u_k$ denoting the $k$-th time derivative of $u$, such that $x$ and $u$ can locally be parameterized by $y$ and its time derivatives. For this parameterization we write
		\vspace{-0.8ex}
		\begin{align*}
			\begin{aligned}
				x&=F_x(y,y_1,\ldots,y_{r-1})\\
				u&=F_u(y,y_1,\ldots,y_r)
			\end{aligned}
		\end{align*}
		\vspace{-2.5ex}
		
		\noindent
		and refer to it as the parameterizing map with respect to the flat output $y$. For a given flat output, $F_x$ and $F_u$ are unique. If the parameterizing map is invertible, \ie $y$ and all the time derivatives of $y$ present in the map can be expressed solely as functions of $x$ and $u$, the system is static feedback linearizable. In this case we call $y$ a linearizing output of the static feedback linearizable system. The static feedback linearization problem has been completely solved, see \cite{JakubczykRespondek:1980,NijmeijervanderSchaft:1990}. However, there do not exist verifiable necessary and sufficient conditions for flatness, except for certain classes of systems, including driftless systems with two inputs \cite{MartinRouchon:1994} and systems which are linearizable by a one-fold prolongation of a suitably chosen control \cite{NicolauRespondek:2017}. Necessary and sufficient conditions for $(x,u)$-flatness of two-input control affine systems with four states can be found in \cite{Pomet:1997}.
		
		It is well known that every flat system can be rendered static feedback linearizable by an endogenous dynamic feedback. If a flat output is known, such a linearizing feedback can be constructed systematically, see \eg \cite{FliessLevineMartinRouchon:1999}. In this contribution we propose an algorithmic test to check whether a two input system is linearizable by an endogenous dynamic feedback with a dimension of at most two. This test furthermore provides a procedure for systematically deriving flat outputs of such systems. The main idea of our algorithmic test is to successively split off or add endogenous dynamic feedbacks to the system in order to obtain a static feedback linearizable system. If the algorithm does not yield a static feedback linearizable system, \ie if the test fails, the considered system is not linearizable by an endogenous dynamic feedback with a dimension of at most two. The linearizing outputs of the eventually obtained static feedback linearizable system are then flat outputs of the original system. In the literature, constructive approaches for deriving flat outputs can be found, see \eg \cite{SchlacherSchoberl:2013}, \cite{SchoberlSchlacher:2014} or \cite{SchlacherSchoberlKolar:2015}. All these procedures are applicable to a broad variety of systems and are not restricted to the two-input case. However, there usually occur degrees of freedom, which may lead to an infinite number of branches. In contrast to that, in every step of our procedure, a certain action has to be applied depending on easily verifiable properties of the system. We will see that systems which are linearizable by an endogenous dynamic feedback with a dimension of at most two are actually linearizable by a special subclass of endogenous dynamic feedbacks, namely prolongations of a suitably chosen control. In \cite{NicolauRespondek:2016}, a complete solution for the flatness problem of two-input systems linearizable via a one-fold prolongation of a suitably chosen control is provided. In \cite{NicolauRespondek:2016-2}, two-input systems linearizable via a two-fold prolongation of a suitably chosen control are considered. The present contribution is greatly influenced by these results. However, in \cite{NicolauRespondek:2016-2} a complete solution for the flatness problem of two-input systems linearizable via a two-fold prolongation is not provided, due to Assumption 2 therein. Our procedure also applies to the cases which are not covered by the results in \cite{NicolauRespondek:2016-2}. In Section \ref{se:examples}, we apply our procedure to three examples, to none of which the results in \cite{NicolauRespondek:2016-2} are applicable. In \cite{GstottnerKolarSchoberl:2020c} we proposed a structurally flat triangular form based on the extended chained form. Our first two examples turn out to be static feedback equivalent to the proposed triangular form, which provides another way for systematically deriving flat outputs for our first two examples. However, the results in \cite{GstottnerKolarSchoberl:2020c} are not applicable to our third example.
		\vspace{-0.8ex}
	\section{NOTATION}
		\vspace{-0.8ex}
		\noindent
		Let $\mathcal{X}$ be an $n$-dimensional smooth manifold, equipped with local coordinates $x^i$, $i=1,\ldots,n$. Its tangent bundle is denoted by $(\mathcal{T}(\mathcal{X}),\tau_\mathcal{X},\mathcal{X})$, for which we have the induced local coordinates $(x^i,\dot{x}^i)$ with respect to the basis $\{\partial_{x^i}\}$. We make use of the Einstein summation convention. By $\partial_xh$ we denote the $m\times n$ Jacobian matrix of $h=(h^1,\ldots,h^m)$ with respect to $x=(x^1,\ldots,x^n)$. The $k$-fold Lie derivative of a function $\varphi$ along a vector field $v$ is denoted by $\Lie_v^k\varphi$. Let $v$ and $w$ be two vector fields. Their Lie bracket is denoted by $[v,w]$. Let furthermore $D_1$ and $D_2$ be two distributions. By $[v,D_1]$ we denote the distribution spanned by the Lie bracket of $v$ with all basis vector fields of $D_1$, and by $[D_1,D_2]$ the distribution spanned by the Lie brackets of all possible pairs of basis vector fields of $D_1$ and $D_2$. The $i$-th derived flag of a distribution $D$ is denoted by $D^{(i)}$ and defined by $D^{(0)}=D$ and $D^{(i+1)}=D^{(i)}+[D^{(i)},D^{(i)}]$ for $i\geq 0$. We denote the Cauchy characteristic distribution of $D$ by $\C{D}$. It is spanned by all vector fields $c$ which belong to $D$ and satisfy $[c,D]\subset D$. We make use of multi-indices, in particular, by $R=(r_1,r_2)$ we denote the unique multi-index associated to a flat output of a system with two inputs, where $r_j$ denotes the order of the highest derivative of $y^j$ needed to parameterize $x$ and $u$ by this flat output, \ie $y_{[R]}=(y^1,y^1_1,\ldots,y^1_{r_1},y^2,y^2_1,\ldots,y^2_{r_2})$. Furthermore, we define $R\pm c=(r_1\pm c,r_2\pm c)$ with an integer $c$, and $\#R=r_1+r_2$.
		\vspace{-1.2ex}
	\section{PRELIMINARIES}\label{se:preliminaries}
		\vspace{-0.7ex}
		\noindent
		%In the following, we state our procedure for deriving flat outputs of systems which are linearizable by an endogenous dynamic feedback with a dimension of at most two. 
		Throughout, we assume all distributions to have locally constant dimension, we consider generic points only. Consider a nonlinear two-input system of the form
		\vspace{-1ex}
		\begin{align}\label{eq:sys}
			\begin{aligned}
				\dot{x}&=f(x,u)\,,
			\end{aligned}
		\end{align}
		\vspace{-3.3ex}
		
		\noindent
		with $\Dim{x}=n$, $\Dim{u}=2$ and $\Rank{\partial_uf}=2$. Let $y^j=\varphi^j(x,u,u_1,\ldots,u_q)$, $j=1,2$ be a flat output of \eqref{eq:sys} and
		\vspace{-2ex}
		\begin{align}\label{eq:param}
			\begin{aligned}
				x&=F_{x}(y_{[R-1]})\\
				u&=F_{u}(y_{[R]})
			\end{aligned}
		\end{align}
		\vspace{-2.2ex}
		
		\noindent
		the parameterizing map with respect to this flat output. The parameterizing map \eqref{eq:param} is a submersion. We denote the difference of the dimensions of the domain and the codomain of \eqref{eq:param} by $d$, \ie $d=\#R+2-(n+2)=\#R-n$. In \cite{NicolauRespondek:2016} and \cite{NicolauRespondek:2016-2}, the number $\#R+2$ is called the differential weight of the flat output. (The differential weight of a flat output with difference $d$ is thus given by $d+n+2$.) The difference $d$ is the minimal dimension of an endogenous dynamic feedback needed to render \eqref{eq:sys} static feedback linearizable such that $y$ forms a linearizing output of the closed loop system. The endogenous feedback can be constructed systematically, see \eg \cite{FliessLevineMartinRouchon:1999}. If $d=0$, the map \eqref{eq:param} degenerates to a diffeomorphism and the system is static feedback linearizable with $y$ being a linearizing output. A flat output $y$ is called a minimal flat output if the difference $d$ is minimal compared to all other possible flat outputs of the system. We define the difference $d$ of a flat system to be the difference of a minimal flat output. The difference $d$ of a system \eqref{eq:sys} therefore measures its distance from static feedback linearizability, \ie $d$ is the minimal possible dimension of an endogenous dynamic feedback needed to render the system static feedback linearizable. For \eqref{eq:sys}, we define the distributions $D_0=\Span{\partial_{u^1},\partial_{u^2}}$ and $D_i=D_{i-1}+[f,D_{i-1}]$, $i\geq 1$ on the state and input manifold $\mathcal{X}\times\mathcal{U}$, where $f=f^i(x,u)\partial_{x^i}$. 
		\begin{theorem}\label{thm:sfl}
			The two-input system \eqref{eq:sys} is linearizable by static feedback if and only if all the distributions $D_i$ are involutive and $\Dim{D_{n-1}}=n+2$.
		\end{theorem}
	
		For a proof of this theorem, we refer to \cite{NijmeijervanderSchaft:1990}. The involutivity of $D_{i-1}$ implies the invariance of $D_i$ with respect to regular input transformations $\bar{u}=\Phi_u(x,u)$. In particular, since $D_0=\Span{\partial_{u^1},\partial_{u^2}}$ is always involutive, $D_1=D_0+[f,D_0]$ is always feedback invariant, also when the system is not static feedback linearizable. A special case of \eqref{eq:sys} are affine input systems (AI-systems)
		\vspace{-1ex}
		\begin{align}\label{eq:sysAI}
			\begin{aligned}
				\dot{x}&=a(x)+b_1(x)u^1+b_2(x)u^2\,.
			\end{aligned}
		\end{align}
		\vspace{-3.8ex}
		\begin{lemma}\label{lem:ai}
			The system \eqref{eq:sys} allows an AI representation \eqref{eq:sysAI} if and only if $D_0\subset\C{D_1}$.
		\end{lemma}
	
		It follows from the proof of this lemma, which can be found in Section \ref{se:proofs}, that the input transformation $\bar{u}^j=f^j(x,u)$, $j=1,2$, which may require a renumbering of the state variables, yields an AI representation if an AI representation is indeed possible. For AI-systems, we can omit $\partial_{u^j}$ in $D_1$, \ie for AI-systems we define $D_1=\Span{b_1,b_2}$, where $b_j=b_j^i(x)\partial_{x^i}$ are vector fields on the state manifold $\mathcal{X}$. A central role in this contribution plays the more general form
		\vspace{-2.5ex}
		\begin{align}\label{eq:sysPAI}
			\begin{aligned}
				\dot{x}&=a(x,u^1)+b(x,u^1)u^2\,,
			\end{aligned}
		\end{align}
		\vspace{-3.7ex}
		
		\noindent
		in which at least the input $u^2$ occurs affine. This form was introduced in \cite{SchlacherSchoberl:2013} and we refer to it as partial affine input form (PAI form)\footnote{In \cite{Kolar:2017} it is shown that the existence of a PAI representation is actually a necessary condition for flatness and this necessary condition is not restricted to systems with two inputs. See also \cite{Sluis:1993}.}.
		\begin{lemma}\label{lem:pai}
			The system \eqref{eq:sys} allows a PAI representation \eqref{eq:sysPAI} if and only if the condition
			\vspace{-1.4ex}
			\begin{align}\label{eq:alg}
				\begin{aligned}
					(\alpha^1)^2\partial_{u^1}^2f+2\alpha^1\alpha^2\partial_{u^1}\partial_{u^2}f+(\alpha^2)^2\partial_{u^2}^2f&\overset{!}{\in} D_1\,,
				\end{aligned}
			\end{align}
			\vspace{-3.2ex}
			
			\noindent
			admits nontrivial solutions $\alpha^j(x,u)$ and at least one of the vector fields $v_c=\alpha^1\partial_{u^1}+\alpha^2\partial_{u^2}$ obtained from a solution of \eqref{eq:alg} meets $v_c\in\C{\Delta_1}$ where $\Delta_1=\Span{\partial_{u^1},\partial_{u^2},[v_c,f]}$.
		\end{lemma}
	
		A proof of this lemma is provided in Section \ref{se:proofs}, where we also discuss a crucial property of the condition \eqref{eq:alg}, namely that for systems which do not allow an AI representation, \ie $D_0\not\subset\C{D_1}$, the condition \eqref{eq:alg} admits at most two independent nontrivial solutions and thus, it yields at most two non-collinear candidates for a vector field $v_c$. To each vector field $v_c$ which actually fulfills $v_c\in\C{\Delta_1}$, there exists a corresponding PAI representation. To derive it, apply the following two steps:
		\begin{enumerate}[leftmargin=*]
			\item Straighten out the distribution $\Span{v_c}$, \ie apply a suitable input transformation $\bar{u}=\Phi_u(x,u)$ such that $\Span{v_c}=\Span{\partial_{\bar{u}^2}}$.
			\item Normalize one equation which explicitly depends on $\bar{u}^2$ (\ie introduce \eg $\tilde{u}^2=\bar{f}^1(x,\bar{u}^1,\bar{u}^2)$, which may require a renumbering of the states of the system).
		\end{enumerate}
		Based on the fact that at most two non-collinear vector fields $v_c$ exist and that to every PAI representation of the system, a vector field $v_c$ belongs, it can be shown that there exist at most two fundamentally different PAI representations of the system, all others are equivalent to one of those two by a transformation of the form
		\vspace{-0.5ex}
		\begin{align}\label{eq:paiPreservingTransformation}
			\begin{aligned}
				\bar{x}&=\Phi_x(x)\\
				\bar{u}^1&=g^1(x,u^1)\\
				\bar{u}^2&=g^2(x,u^1)+m(x,u^1)u^2\,.
			\end{aligned}
		\end{align}
		\section{FINITE ALGORITHMIC TEST}\noindent
			The main idea of our test for flatness with $d\leq 2$ is to successively split off or add endogenous dynamic feedbacks to the system, in order to eventually obtain a static feedback linearizable system. The linearizing outputs of the eventually obtained static feedback linearizable system are then flat outputs with $d\leq 2$ of the original system. As we will see in the next section, for systems with $d\leq 2$, there always applies one of the following cases.
			\begin{enumerate}[leftmargin=*]
				\item The system allows an AI representation \eqref{eq:sysAI} with an involutive input distribution $D_1=\Span{b_1,b_2}$. In this case, straighten out $D_1$ by a suitable state transformation $\bar{x}=\Phi_x(x)$, in order to obtain a decomposition of \eqref{eq:sysAI} into the form
				\vspace{-0.9ex}
				\begin{align}\label{eq:decompCase1}
					\begin{aligned}
						\Sigma_1:&&\dot{\bar{x}}_1^{i_1}&=\bar{f}_1^{i_1}(\bar{x}_1,\bar{x}_2)\,,&i_1&=1,\ldots,n-2\\
						\Sigma_2:&&\dot{\bar{x}}_2^{i_2}&=\bar{f}_2^{i_1}(\bar{x}_1,\bar{x}_2,u)\,,&i_2&=1,2\,.
					\end{aligned}
				\end{align}
				\vspace{-2.4ex}
				
				\noindent
				Continue the analysis with the subsystem $\Sigma_1$ with the state $\bar{x}_1$ and the input $\bar{x}_2$.\label{case:1}
				\item The system allows an AI representation \eqref{eq:sysAI} but the input distribution $D_1=\Span{b_1,b_2}$ is non-involutive. Two subcases are possible. If $\Dim{\overline{D}_1}=3$, construct a nontrivial linear combination $b_c=\alpha^1b_1+\alpha^2b_2$ such that $[a,b_c]\in\overline{D}_1$. If $\Dim{\overline{D}_1}=4$, construct a nontrivial linear combination $b_c=\alpha^1b_1+\alpha^2b_2$ such that $b_c\in\C{D_1^{(1)}}$. In either subcase, introduce the new input $\bar{u}^1=\alpha^2u^1-\alpha^1u^2$ and one-fold prolong it, \ie add the equation $\dot{\bar{u}}^1=\bar{u}^1_1$. Continue the analysis with the prolonged system
				\vspace{-3.2ex}
				\begin{align}\label{eq:prolongationCase2}
					&\begin{aligned}
						\dot{x}&=a(x)+\bar{b}_1(x)\bar{u}^1+\bar{b}_2(x)\bar{u}^2\,,
					\end{aligned}&\begin{aligned}
						\dot{\bar{u}}^1&=\bar{u}^1_1
					\end{aligned}
				\end{align}
				\vspace{-3.5ex}
				
				\noindent
				with the state $(x,\bar{u}^1)$ and the input $(\bar{u}^1_1,\bar{u}^2)$.\label{case:2}
				\item The system only allows a PAI representation \eqref{eq:sysPAI}. Transform the system into PAI form and one-fold prolong the non-affine occurring input $\bar{u}^1$, \ie add the equation $\dot{\bar{u}}^1=\bar{u}^1_1$.
				Continue the analysis with the prolonged system
				\vspace{-3.2ex}
				\begin{align}\label{eq:prolongationCase3}
					&\begin{aligned}
						\dot{x}&=a(x,\bar{u}^1)+b(x,\bar{u}^1)\bar{u}^2\,,
					\end{aligned}&\begin{aligned}
						\dot{\bar{u}}^1&=\bar{u}^1_1
					\end{aligned}
				\end{align}
				\vspace{-3.5ex}
				
				\noindent
				with the state $(x,\bar{u}^1)$ and the input $(\bar{u}^1_1,\bar{u}^2)$. If two non-equivalent PAI representation exist, we have to continue with both of them (branching point).\label{case:3}
			\end{enumerate}
			The procedure succeeds if a static feedback linearizable system is obtained. %and fails if in some step none of the cases applies, \ie if the system allows an AI representation but $D_1$ is non-involutive and no linear combination $b_c$ needed in the subcases of case \ref{case:2} exists, or if the system does not even allow a PAI representation. 
			In the next section, we explain in detail the effects of the individual steps which are applied depending on which case applies. %Based on that, we will show that for systems with $d\leq 2$, the procedure necessarily succeeds and conversely, that if the procedure succeeds, the system necessarily has a difference of $d\leq 2$. 
			As we will see, the subsystem $\Sigma_1$ in \eqref{eq:decompCase1} which is obtained if case \ref{case:1} applies, has the same difference $d$ as the complete system. Such a step therefore only reduces the dimension of the state, but not the distance from static feedback linearizability. In contrast to that, if for a system with $d\leq 2$ the cases 2 or 3 apply, we obtain a prolonged system \eqref{eq:prolongationCase2} or \eqref{eq:prolongationCase3} with a difference of $d-1$ compared to the non-prolonged system. If case \ref{case:3} applies, we may encounter a branching point since there may exist two non-equivalent PAI representations which both have to be considered. Nevertheless, at least one of the AI-systems obtained by prolonging the non-affine input has a difference of $d-1$ compared to the non-prolonged system. Therefore, applying the procedure to a system with $d\leq 2$, after at most two steps which actually decrease the difference \mbox{(cases \ref{case:2} and \ref{case:3})}, we must obtain a static feedback linearizable system. Otherwise, the original system has a difference of $d\geq 3$ or is not flat. %Although for system with $d\geq 3$ the procedure may succeed if we do not stop prematurely, we cannot guarantee that the steps \ref{case:2} and \ref{case:3} indeed yield a system with a lower difference than the original system. 
			\vspace{-1.3ex}
		\subsection{Necessity and sufficiency for systems with $d\leq 2$}\noindent
			\vspace{-3.2ex}
			
			\noindent
			In the following, we explain why applying the procedure to a system with a difference of $d\leq 2$ necessarily yields a static feedback linearizable system, and why the linearizing outputs of the obtained static feedback linearizable system are flat outputs with $d\leq 2$ of the original system. We have the following three results, which are proven in Section \ref{se:proofs}.
			\begin{lemma}\label{lem:case1}
				Every flat output with difference $d$ of the complete system \eqref{eq:decompCase1} is also a flat output with the same difference $d$ for the subsystem $\Sigma_1$ in \eqref{eq:decompCase1}. Conversely, every flat output of $\Sigma_1$ is also a flat output of the complete system and the differences again coincide. 
			\end{lemma}
			\begin{lemma}\label{lem:case2}
				For an AI-system with a difference of $d\leq 2$ which has a non-involutive input distributions $D_1=\Span{b_1,b_2}$, we either have $\Dim{\overline{D}_1}=3$ and there exists a non-trivial vector field $b_c=\alpha^1b_1+\alpha^2b_2$ satisfying $[a,b_c]\in\overline{D}_1$, or we have $\Dim{\overline{D}_1}=4$ and there exists a non-trivial vector field $b_c=\alpha^1b_1+\alpha^2b_2$ satisfying $b_c\in\C{D_1^{(1)}}$. Furthermore, the prolonged system \eqref{eq:prolongationCase2} has a difference of $d-1$ (\ie it is either static feedback linearizable or has a difference of $1$). The linearizing outputs or flat outputs with difference $1$ of the prolonged system are flat outputs with difference $1$ or $2$ of the original system.
			\end{lemma}
			\begin{lemma}\label{lem:case3}
				A system with $d\leq 2$ which does not allow an AI representation allows a PAI representation. There exist at most two non-equivalent PAI representations and at least one of the prolonged systems \eqref{eq:prolongationCase3} obtained from these possibly two PAI representations has a difference of $d-1$ (\ie it is either static feedback linearizable or has a difference of $1$). The linearizing outputs or flat outputs with difference $1$ of the prolonged system are flat outputs with difference $1$ or $2$ of the original system.
			\end{lemma}
		
			Based on these results, we first explain why the procedure necessarily succeeds when it is applied to a system with a difference of $d\leq 2$. Assume we have a system \eqref{eq:sys} with $d\leq 2$. If the system allows an AI representation, its input distribution $D_1=\Span{b_1,b_2}$ can either be involutive or not. If it is involutive, case \ref{case:1} applies, which according to Lemma \ref{lem:case1} yields a system with the same difference $d$ as the original system. We can exclude the case $\Dim{D_1}\leq 1$, since then we effectively have a single input system which is either static feedback linearizable or not flat at all, see \cite{CharletLevineMarino:1991}. Furthermore, we can exclude the case that the subsystem $\Sigma_1$ in \eqref{eq:decompCase1} has redundant inputs\footnote{If $\Sigma_1$ has redundant inputs, \ie if $\Sigma_1$ effectively reduces to a single input system, it can either be static feedback linearizable or not flat at all. However, if $\Sigma_1$ is static feedback linearizable, then also the complete system is static feedback linearizable and the procedure already succeeded before this decomposition. On the other hand, if $\Sigma_1$ is not flat then due to Lemma \ref{lem:case1}, also \eqref{eq:decompCase1} is not flat. The system \eqref{eq:decompCase1} need not be the original system, it may be the outcome of previously applied decompositions or prolongations according to the cases 1 to 3. Nevertheless, due to the Lemmas \ref{lem:case1} to \ref{lem:case3}, if $\Sigma_1$ is not flat, also the original system cannot be flat with a difference of $d\leq2$. In fact, it even follows that the original system cannot be flat at all.}. If the procedure succeeds, it necessarily terminates with case \ref{case:2} or \ref{case:3} since case \ref{case:1} always yields a subsystem $\Sigma_1$ with the same difference as the complete system \eqref{eq:decompCase1}.
			
			If the system allows an AI representation with a non-involutive input distribution $D_1$, case \ref{case:2} applies and due to Lemma \ref{lem:case2}, the prolonged system \eqref{eq:prolongationCase2} has a difference of $d\leq 1$. 
			
			If the system does not allow an AI representation, according to Lemma \ref{lem:case3} it necessarily allows a PAI representation and thus, case \ref{case:3} applies. Since there exist at most two non-equivalent PAI representations and in case \ref{case:3} we require to consider both of them, one of the possibly two prolonged systems \eqref{eq:prolongationCase3} has a difference of $d\leq 1$.
			
			Each step therefore either yields a system with a lower state dimension or a lower difference (for at least one of the possibly two branches if case 3 applies). Since by assumption the original system has a difference of $d\leq 2$, after at most two steps which actually yield a system with a lower difference (cases \ref{case:2} and \ref{case:3}) we necessarily have a static feedback linearizable system.
			
			On the other hand, due to Lemma \ref{lem:case1} to \ref{lem:case3} it follows that if the procedure succeeds, \ie if after at most two steps which actually yield a system with a lower difference (cases \ref{case:2} and \ref{case:3}) we obtain a static feedback linearizable system, the linearizing outputs of this static feedback linearizable system are flat outputs with $d=1$ or $d=2$ of the original system.
	\section{EXAMPLES}\label{se:examples}
		\noindent
		As already mentioned in the introduction, our first two examples can also be handled with the structurally flat triangular form which we proposed in \cite{GstottnerKolarSchoberl:2020c}. However, neither \cite{GstottnerKolarSchoberl:2020c} nor \cite{NicolauRespondek:2016-2}, nor both of them together completely cover the class of systems addressed in this paper. The results in \cite{GstottnerKolarSchoberl:2020c} and \cite{NicolauRespondek:2016-2} do not apply to our third example.
		\subsection{Planar VTOL aircraft}
			\noindent
			Consider the planar VTOL aircraft, also treated in \eg \cite{FliessLevineMartinRouchon:1999}, \cite{SchoberlRiegerSchlacher:2010}, \cite{SchoberlSchlacher:2011} or \cite{GstottnerKolarSchoberl:2020c} and given by
			\begin{align}\label{eq:vtol}
				\begin{aligned}
					\dot{x}&=v_x&&&\dot{v}_x&=\epsilon\cos(\theta)u^2-\sin(\theta)u^1\\
					\dot{z}&=v_z&&&\dot{v}_z&=\cos(\theta)u^1+\epsilon\sin(\theta)u^2-1\\
					\dot{\theta}&=\omega&&&\dot{\omega}&=u^2\,.
				\end{aligned}
			\end{align}
			The input vector fields of this system are given by $b_1=-\sin(\theta)\partial_{v_x}+\cos(\theta)\partial_{v_z}$ and $b_2=\epsilon\cos(\theta)\partial_{v_x}+\epsilon\sin(\theta)\partial_{v_z}+\partial_\omega$, which span the involutive input distribution $D_1=\Span{b_1,b_2}$. Therefore, case \ref{case:1} of the procedure applies. In order to decompose the system as described in case \ref{case:1}, we straighten out $D_1$ by applying the state transformation $\bar{v}_x=\cos(\theta)v_x+\sin(\theta)v_z-\epsilon\omega$, with the rest of the coordinates left unchanged. This results in the decomposition
			\begin{align}\label{eq:vtolDecomp}
				\begin{aligned}
					&\Sigma_1:\quad\begin{aligned}
						\dot{x}&=\tfrac{1}{\cos(\theta)}(\bar{v}_x-\sin(\theta)v_z+\epsilon\omega)\\
						\dot{z}&=v_z\\
						\dot{\theta}&=\omega\\
						\dot{\bar{v}}_x&=\tfrac{\omega}{\cos(\theta)}(v_z-\sin(\theta)(\bar{v}_x+\epsilon\omega))-\sin(\theta)\\
					\end{aligned}\\[0.5ex]
					&\Sigma_2:\quad\begin{aligned}
						\dot{v}_z&=\cos(\theta)u^1+\epsilon\sin(\theta)u^2-1\\
						\dot{\omega}&=u^2\,.
					\end{aligned}
				\end{aligned}
			\end{align}
			We proceed with the subsystem $\Sigma_1$, for which $v_z$ and $\omega$ act as inputs. This system does not allow an AI representation, so we are in case \ref{case:3} of the procedure. In order to derive a PAI representation for $\Sigma_1$, we first solve \eqref{eq:alg}, which yields $\alpha^2(\alpha^1-\alpha^2\epsilon\sin(\theta))\overset{!}{=}0$, and has the two independent solutions $\alpha^1=\lambda$, $\alpha^2=0$ and $\alpha^1=\lambda\epsilon\sin(\theta)$, $\alpha^2=\lambda$, with an arbitrary function $\lambda\neq 0$ (we can choose \eg $\lambda=1$). Both vector fields $v_{c,1}=\partial_{v_z}$ and $v_{c,2}=\epsilon\sin(\theta)\partial_{v_z}+\partial_\omega$ obtained from these solutions fulfill $v_{c,i}\in\C{\Span{\partial_{v_z},\partial_\omega,[v_{c,i},f]}}$, \ie $\Sigma_1$ allows two non-equivalent PAI representations. The vector field $v_{c,1}=\partial_{v_z}$ is already straightened out and $\Sigma_1$ in \eqref{eq:vtolDecomp} is actually already in PAI form with the input $v_z$ occurring affine. However, the AI-system obtained by one-fold prolonging the non-affine input $\omega$
%			, \ie
%			\begin{align}\label{eq:vtolAI1}
%				\begin{aligned}
%					\dot{x}&=\tfrac{1}{\cos(\theta)}(\bar{v}_x-\sin(\theta)v_z+\epsilon\omega)\\
%					\dot{z}&=v_z\\
%					\dot{\theta}&=\omega\\
%					\dot{\bar{v}}_x&=\tfrac{\omega}{\cos(\theta)}(v_z-\sin(\theta)(\bar{v}_x-\epsilon\omega))-\sin(\theta)\\
%					\dot{\omega}&=\omega_1\,,
%				\end{aligned}
%			\end{align}
			cannot have a difference of $d=1$, since none of the cases of the procedure applies to this AI-system.
%			(For \eqref{eq:vtolAI1}, the distribution $D_1=\Span{b_1,b_2}$, spanned by its input vector fields $b_1=-\tfrac{\sin(\theta)}{\cos(\theta)}\partial_x+\partial_z+\tfrac{\omega}{\cos(\theta)}\partial_{\bar{v}_x}$ and $b_2=\partial_\omega$ is non-involutive and $\Dim{\overline{D}_1}=3$, but there does not exist a non-trivial linear combination $b_c=\alpha^1b_1+\alpha^2b_2$ which satisfies $[a,b_c]\in \overline{D}_1$, where $a=\tfrac{v_x+\epsilon\omega}{\cos(\theta)}\partial_x+\omega\partial_\theta+\sin(\theta)(\tfrac{\omega}{\cos(\theta)}(\epsilon\omega-\bar{v}_x)-1)\partial_{\bar{v}_x}$is the drift of \eqref{eq:vtolAI1}.)
			Straightening out $v_{c,2}=\epsilon\sin(\theta)\partial_{v_z}+\partial_\omega$, by the input transformation $\bar{v}_z=v_z-\epsilon\sin(\theta)\omega$, results in $v_{c,2}=\partial_{\omega}$ and the PAI representation\footnote{Note that in these coordinates in fact both vector fields $v_{c,1}$ and $v_{c,2}$ are straightened out simultaneously and that in \eqref{eq:vtolPAI2} actually either input could be interpreted as the affine entering input.}
			\vspace{-0.8ex}
			\begin{align}\label{eq:vtolPAI2}
				\begin{aligned}
					\dot{x}&=\tfrac{1}{\cos(\theta)}(\bar{v}_x-\sin(\theta)\bar{v}_z)+\epsilon\cos(\theta)\omega\\
					\dot{z}&=\bar{v}_z+\epsilon\sin(\theta)\omega\\
					\dot{\theta}&=\omega\\
					\dot{\bar{v}}_x&=\tfrac{\omega}{\cos(\theta)}(\bar{v}_z-\sin(\theta)\bar{v}_x)-\sin(\theta)\,.
				\end{aligned}
			\end{align}
			\vspace{-2.1ex}
			
			\noindent
			For the AI-system obtained by one-fold prolonging $\bar{v}_z$, \ie adding the equation $\dot{\bar{v}}_z=\bar{v}_{z,1}$, we have $\Dim{\overline{D}}=3$, where $D_1=\Span{b_1,b_2}$ is the input distribution spanned by its input vector fields $b_1=\partial_{\bar{v}_z}$ and $b_2=\epsilon\cos(\theta)\partial_x+\epsilon\sin(\theta)\partial_z+\partial_\theta+\tfrac{\bar{v}_z-\sin(\theta)\bar{v}_x}{\cos(\theta)}\partial_{\bar{v}_x}$. The drift of this AI-system is given by $a=\tfrac{1}{\cos(\theta)}(\bar{v}_x-\sin(\theta)\bar{v}_z)\partial_x+\bar{v}_z\partial_z-\sin(\theta)\partial_{\bar{v}_x}$ and there indeed exists a linear combination $b_c=\alpha^1b_1+\alpha^2b_2$ of the input vector fields which satisfies $[a,b_c]\in\overline{D}_1$, namely $b_c=\lambda b_2$ with an arbitrary function $\lambda\neq 0$ (\eg $\lambda=1$ and thus $\alpha^1=0$ and $\alpha^2=1$). So the conditions for case \ref{case:2} are met and we prolong the input $\tilde{u}^1=\alpha^2\bar{v}_{z,1}-\alpha^1\omega=\bar{v}_{z,1}$, \ie we simply have to prolong the input $\bar{v}_{z,1}$, no actual input transformation is required in this case, and we obtain the static feedback linearizable system
			\vspace{-0.8ex}
			\begin{align*}%\label{eq:vtolSfl}
				&\begin{aligned}
					\dot{x}&=\tfrac{1}{\cos(\theta)}(\bar{v}_x-\sin(\theta)\bar{v}_z)+\epsilon\cos(\theta)\omega\\
					\dot{z}&=\bar{v}_z+\epsilon\sin(\theta)\omega\\
					\dot{\theta}&=\omega\\
					\dot{\bar{v}}_x&=\tfrac{\omega}{\cos(\theta)}(\bar{v}_z-\sin(\theta)\bar{v}_x)-\sin(\theta)
				\end{aligned}&\begin{aligned}
					\dot{\bar{v}}_z&=\bar{v}_{z,1}\\
					\dot{\bar{v}}_{z,1}&=\bar{v}_{z,2}\,,
				\end{aligned}
			\end{align*}
			\vspace{-2.1ex}
			
			\noindent
			with the input $(\omega,\bar{v}_{z,2})$. We have $y=(x-\epsilon\sin(\theta),z+\epsilon\cos(\theta))$ as a possible linearizing output of this system, which in turn is a flat output with a difference of $d=2$ of the original system \eqref{eq:vtol}.
			\vspace{-1ex}
		\subsection{Academic example I}
			\vspace{-0.5ex}
			\noindent
			Consider the system
			\vspace{-0.9ex}
			\begin{align}\label{eq:sin}
				\begin{aligned}
					\dot{x}^1&=u^1\,,&\dot{x}^2&=u^2\,,&\dot{x}^3&=\sin(\tfrac{u^1}{u^2})\,,
				\end{aligned}
			\end{align}
			\vspace{-3.6ex}
			
			\noindent
			also considered in \cite{Levine:2009}, \cite{Schoberl:2014} and \cite{GstottnerKolarSchoberl:2020c}. This system does not allow an AI representation, so we are in case \ref{case:3}. In order to derive a PAI representation, we solve \eqref{eq:alg}, which for this system yields
			\begin{align*}
				\begin{aligned}
					&(\alpha^1)^2\sin(\tfrac{u^1}{u^2})(u^2)^2+2\alpha^1\alpha^2(\cos(\tfrac{u^1}{u^2})u^2-\sin(\tfrac{u^1}{u^2})u^1)u^2+\\
					&\hspace{6em}(\alpha^2)^2(\sin(\tfrac{u^1}{u^2})u^1-2\cos(\tfrac{u^1}{u^2})u^2)u^1\overset{!}{=}0
				\end{aligned}
			\end{align*}
			and has the two independent non-trivial solutions $\alpha^1=\lambda u^1$, $\alpha^2=\lambda u^2$ and $\alpha^1=\lambda(u^1\tan(\tfrac{u^1}{u^2})-2u^2)$, $\alpha^2=\lambda u^2\tan(\tfrac{u^1}{u^2})$ with an arbitrary function $\lambda\neq 0$ (\eg $\lambda=1$). Only the vector field $v_c=u^1\partial_{u^1}+u^2\partial_{u^2}$, which is obtained from the first solution, meets $v_c\in\C{\Span{\partial_{u^1},\partial_{u^2},[v_c,f]}}$. To derive the corresponding PAI representation, we straighten out $\Span{v_c}$ by the input transformation $\bar{u}^1=\tfrac{u^1}{u^2}$, $\bar{u}^2=u^2$, resulting in
			\vspace{-0.5ex}
			\begin{align*}
				\begin{aligned}
					\dot{x}^1&=\bar{u}^1\bar{u}^2\,,&\dot{x}^2&=\bar{u}^2\,,&\dot{x}^3&=\sin(\bar{u}^1)
				\end{aligned}
			\end{align*}
			\vspace{-3.3ex}
			
			\noindent
			and proceed with the AI-system obtained by one-fold prolonging $\bar{u}^1$, \ie adding the equation $\dot{\bar{u}}^1=\bar{u}^1_1$. Its input vector fields are given by $b_1=\partial_{\bar{u}^1}$, $b_2=\bar{u}^1\partial_{x^1}+\partial_{x^2}$ and its drift reads $a=\sin(\bar{u}^1)\partial_{x^3}$. We have $\Dim{\overline{D}_1}=3$ and there exists a linear combination $b_c=\alpha^1b_1+\alpha^2b_2$ which meets $[a,b_c]\in\overline{D}_1$, namely $b_c=\lambda b_2$, with an arbitrary function $\lambda\neq 0$ (\eg $\lambda=1$ and thus $\alpha^1=0$ and $\alpha^2=1$). So the conditions for case \ref{case:2} are met and we prolong the input $\tilde{u}^1=\alpha^2\bar{u}^1_1-\alpha^1\bar{u}^2=\bar{u}^1_1$, \ie we simply have to prolong the input $\bar{u}^1_1$, no actual input transformation is required, and we obtain the static feedback linearizable system
			\begin{align}\label{eq:sinSfl}
				&\begin{aligned}
					\dot{x}^1&=\bar{u}^1\bar{u}^2\\
					\dot{x}^2&=\bar{u}^2\\
					\dot{x}^3&=\sin(\bar{u}^1)\\
				\end{aligned}&
				\begin{aligned}
					\dot{\bar{u}}^1&=\bar{u}^1_1\\
					\dot{\bar{u}}^1_1&=\bar{u}^1_2
				\end{aligned}
			\end{align}
			with the input $(\bar{u}^1_2,\bar{u}^2)$. We have $y=(x^3,x^1-x^2\bar{u}^1)=(x^3,x^1-x^2\tfrac{u^1}{u^2})$ as a possible linearizing output of \eqref{eq:sinSfl}, which in turn is a flat output with a difference of $d=2$ of the original system \eqref{eq:sin}.
			\vspace{-1ex}
		\subsection{Academic example II}
			\vspace{-0.5ex}
			\noindent
			As a final example, consider the system
			\vspace{-0.8ex}
			\begin{align*}
				&\begin{aligned}
					\dot{x}^1&=\arcsin(\tfrac{u^1+u^2}{x^2})-x^4\\
					\dot{x}^2&=x^4
				\end{aligned}&\begin{aligned}
					\dot{x}^3&=u^1\\
					\dot{x}^4&=u^2\,.
				\end{aligned}
			\end{align*}
			\vspace{-2.6ex}
			
			\noindent
			With the algorithmic test proposed in this paper, it can be shown that this system has a difference of $d=2$ and $y=(x^1+x^2,x^3+x^4)$ follows as a corresponding flat output with $d=2$. In the first step case \ref{case:3} applies, in the second step case \ref{case:1} applies, and then case \ref{case:3} applies again.

	\section{CONCLUSIONS}
		\vspace{-0.5ex}
		\noindent
		We have proposed a finite algorithmic test for linearizability by an endogenous dynamic feedback with a dimension of at most two, together with a procedure for systematically deriving flat outputs for such systems. The conditions are in principle verifiable, however, applying the test requires straightening out involutive distributions, which from a computational point of view is unfavorable. Further research will be devoted to deriving necessary and sufficient conditions based on the proposed test which overcome this drawback.
		
%		In \cite{Kolar:2017}, it is shown that certain bounds on the orders $r_j$ of the time derivatives of $y$ in the flat parameterization hold, from which it follows that every $x$-flat two-input system with four states and configuration flat system with six states belongs to the class of systems addressed in this paper.
	\addtolength{\textheight}{-0cm}   % This command serves to balance the column lengths
	                                  % on the last page of the document manually. It shortens
	                                  % the textheight of the last page by a suitable amount.
	                                  % This command does not take effect until the next page
	                                  % so it should come on the page before the last. Make
	                                  % sure that you do not shorten the textheight too much.
	
	%%%%%%%%%%%%%%%%%%%%%%%%%%%%%%%%%%%%%%%%%%%%%%%%%%%%%%%%%%%%%%%%%%%%%%%%%%%%%%%%

	%%%%%%%%%%%%%%%%%%%%%%%%%%%%%%%%%%%%%%%%%%%%%%%%%%%%%%%%%%%%%%%%%%%%%%%%%%%%%%%%

	%%%%%%%%%%%%%%%%%%%%%%%%%%%%%%%%%%%%%%%%%%%%%%%%%%%%%%%%%%%%%%%%%%%%%%%%%%%%%%%%
	\vspace{-0.7ex}
	\section{PROOFS}\label{se:proofs}
		\vspace{-0.9ex}
		\subsection*{Proof of Lemma \ref{lem:ai}}\noindent
			\textit{Necessity.} For an AI-system \eqref{eq:sysAI}, we have $D_1=\Span{\partial_{u^1},\partial_{u^2},b_1,b_2}$ and since the vector fields $b_j$ do not depend on $u^1$ and $u^2$, we obviously have $D_0=\Span{\partial_{u^1},\partial_{u^2}}\subset\C{D_1}$.
			
			\noindent			
			\textit{Sufficiency.} Given a system $\dot{x}=f(x,u)$ with $\Rank{\partial_uf}=2$, we can always apply an input transformation $\bar{u}^j=f^j(x,u)$, $j=1,2$ (this transformation may require a renumbering of the sate variables) in order to obtain
			\begin{align}\label{eq:normalized}
				\begin{aligned}
					\dot{x}^j&=\bar{u}^j\,,&j&=1,2\\
					\dot{x}^i&=\bar{f}^i(x,\bar{u})\,,&i&=3,\ldots,n\,.
				\end{aligned}
			\end{align}
			Calculating $D_1$ in these coordinates yields $D_1=\Span{\partial_{\bar{u}^1},\partial_{\bar{u}^2},v_1,v_2}$ with $v_1=\partial_{x^1}+\partial_{\bar{u}^1}\bar{f}^i\partial_{x^i}$, $i=3,\ldots,n$ and $v_2=\partial_{x^2}+\partial_{\bar{u}^2}\bar{f}^i\partial_{x^i}$, $i=3,\ldots,n$. The condition $D_0\subset\C{D_1}$ implies $[\partial_{\bar{u}^j},v_k]\in D_1$, $j,k\in\{1,2\}$, which can only hold if $[\partial_{\bar{u}^j},v_k]=0$, implying that \eqref{eq:normalized} is actually in AI form.\hfill\QED
		\subsection*{Proof of Lemma \ref{lem:pai}}\noindent
			For a system in PAI form \eqref{eq:sysPAI}, we have $D_1=\Span{\partial_{u^1},\partial_{u^2},(\partial_{u^1}a^i+u^2\partial_{u^1}b^i)\partial_{x^i},b^i\partial_{x^i}}$. Consider the subdistribution $\Delta_1=\Span{\partial_{u^1},\partial_{u^2},b^i\partial_{x^i}}$ of $D_1$. We have $\partial_{u^2}\in\C{\Delta_1}$ and $\Delta_1=D_0+\Span{[\partial_{u^2},f]}$. Therefore, if a two input system \eqref{eq:sys} allows a PAI representation, then there also must exist a vector field $v_c\in D_0$, \ie $v_c=\alpha^1\partial_{u^1}+\alpha^2\partial_{u^2}$ with functions $\alpha^j=\alpha^j(x,u)$, such that with $\Delta_1=D_0+\Span{[v_c,f]}$, we have $v_c\in\C{\Delta_1}$. Since $v_c\in D_0$ and $D_0$ is involutive, this is equivalent to $[v_c,[v_c,f]]\in\Delta_1$. Since we have $\Delta_1\subset D_1$, this also implies $[v_c,[v_c,f]]\in D_1$. Inserting $v_c=\alpha^1\partial_{u^1}+\alpha^2\partial_{u^2}$ into the latter relation yields \eqref{eq:alg}. On the other hand, given a vector field $v_c=\alpha^1\partial_{u^1}+\alpha^2\partial_{u^2}$ which meets $[v_c,[v_c,f]]\in\Delta_1$, a PAI representation can be derived following the two steps in Section \ref{se:preliminaries}. To show that this procedure indeed leads to a PAI representation, evaluate the condition $[v_c,[v_c,f]]\in\Delta_1$ after applying the two steps of the procedure in Section \ref{se:preliminaries}.  
			
			We do not show in detail here why the condition \eqref{eq:alg} admits at most two independent nontrivial solutions for systems which do not allow an AI representation, \ie $D_0\not\subset\C{D_1}$. The proof is based on the fact that $D_0\not\subset\C{D_1}$ implies that at least one of the vector fields $\partial_{u^1}^2f$, $\partial_{u^1}\partial_{u^2}f$ or $\partial_{u^2}^2f$ in \eqref{eq:alg} is not contained in $D_1$. If non of them is contained in $D_1$, it follows that \eqref{eq:alg} admits no non-trivial solution. If one of them is$\Mod D_1$ linearly dependent of the other two, \eqref{eq:alg} admits at most one independent non-trivial solution. If two of them are$\Mod D_1$ linearly dependent of one of them, \eqref{eq:alg} admits at most two independent non-trivial solutions. %Possibly existing complex solutions of \eqref{eq:alg} are not relevant, only real ones are counted as solutions.
			\hfill\QED
			\vspace{-0.5ex}
		\subsection*{Proof of Lemma \ref{lem:case1}}
			\vspace{-0.5ex}
			\noindent
			Given a flat output $y$ of the complete system \eqref{eq:decompCase1}, we have a certain difference $d=\#R-n$ (which by assumption is at most two) of the dimensions of the domain and the codomain of the corresponding parameterizing map
			\vspace{-0.4ex}
			\begin{subequations}
				\begin{align}
				(x_1,x_2)&=F_x(y_{[R-1]})\label{eq:xpar}\\		
				u&=F_u(y_{[R]})\,.
				\end{align}
			\end{subequations}
			The parameterizing map of the subsystem $\Sigma_1$ in \eqref{eq:decompCase1} only consists of \eqref{eq:xpar}, \ie $y$ is also a flat output of the subsystem $\Sigma_1$%\footnote{If $y=\varphi(\bar{x},u,u_1,u_{1,1},\ldots,u_2,u_{2,1},\ldots)$ actually depends on $u$ or time derivatives $u_{\alpha}$, $\alpha \geq 1$ of $u$, those have to be expressed in terms of $\bar{x}_1$, $\bar{x}_2$ and time derivatives of $\bar{x}_2$ in order to obtain a flat output which only depends on the state $\bar{x}_1$, input $\bar{x}_2$ and time derivatives of the input $\bar{x}_2$ of the subsystem $\Sigma_1$. However, it turns out that an AI-system with $d\leq 2$ is actually always $x$-flat, see Lemma \ref{lem:d2xu} in Section \ref{se:proofs}.}
			. For the parameterizing map \eqref{eq:xpar} of $\Sigma_1$ we have the same difference of the dimensions of the domain and the codomain, \ie the same value for $d$. Indeed, in \eqref{eq:xpar}, we have a domain of dimension $\#(R-1)+2=\#R$ and a codomain of dimension $n$ and thus again a difference of $d=\#R-n$. Similarly, it can be shown that every flat output of $\Sigma_1$ with a certain difference $d$ is also a flat output for the complete system with the same difference $d$.\hfill\QED
			\vspace{1ex}
			
			The proofs of Lemma \ref{lem:case2} and \ref{lem:case3} are based on the following Theorem, which is proven at the end of this section.
			\begin{theorem}\label{thm:d2lin}
				A system \eqref{eq:sys} with $d\leq 2$ can be rendered static feedback linearizable by $d$-fold prolonging a suitably chosen (new) input after a suitable input transformation has been applied. If the system under consideration is an AI-system, there always exists an affine input transformation which generates the required (new) input. %vielleicht passt hier "(new) control afte a suitable static feedback has been applied" besser
			\end{theorem}
%			\begin{lemma}\label{lem:minimal}
%				Let $y$ be a minimal flat output of the system \eqref{eq:sys} with a certain (minimal) difference $d$. Assume that there exists an input transformation $\bar{u}=\Phi_u(x,u)$ such that the parameterization of the new input $\bar{u}^1$ involves derivatives of $y$ up to order $R-1$ only, \ie $\bar{u}^1=F_u^1(y_{[R-1]})$. Then, $y$ is also a minimal flat output of the prolonged system
%				\begin{align}\label{eq:minimalProlonged}
%					\begin{aligned}
%						\dot{x}&=\bar{f}(x,\bar{u}^1,\bar{u}^2)\\
%						\dot{\bar{u}}^1&=\bar{u}^1_1
%					\end{aligned}
%				\end{align}
%				with the state $x_p=(x,\bar{u}^1)$ and the input $(\bar{u}^1_1,\bar{u}^2)$, which has a difference of $d-1$. Conversely, every flat output of \eqref{eq:minimalProlonged} with $d-1$ (\ie every minimal flat output of \eqref{eq:minimalProlonged}) is also a minimal flat output of the original system.
%			\end{lemma}
		\vspace{-0.5ex}
		\subsection*{Proof of Lemma \ref{lem:case2}}\noindent
			We prove the cases $d=1$ and $d=2$ separately. Let us first assume that the AI-system has a difference of $d=1$. According to Theorem \ref{thm:d2lin}, it becomes static feedback linearizable after applying a suitable affine input transformation $\bar{u}^j=g(x)^j+m^j_k(x)u^k$ and subsequently one-fold prolonging the new input $\bar{u}^1$, \ie the prolonged system
			\vspace{-0.5ex}
			\begin{align*}
				&\begin{aligned}
					\dot{x}&=\bar{a}(x)+\bar{b}_1(x)\bar{u}^1+\bar{b}_2(x)\bar{u}^2\,,
				\end{aligned}&\begin{aligned}
					\dot{\bar{u}}^1&=\bar{u}^1_1
				\end{aligned}
			\end{align*}
			\vspace{-2.9ex}
			
			\noindent
			with the state $(x,\bar{u}^1)$ and the input $(\bar{u}^1_1,\bar{u}^2)$ is static feedback linearizable. Thus, the distributions $\Delta_1=\Span{\partial_{\bar{u}^1},\bar{b}_2}$, $\Delta_2=\Span{\partial_{\bar{u}^1},\bar{b}_1,\bar{b}_2,[\bar{a},\bar{b}_2]+\bar{u}^1[\bar{b}_1,\bar{b}_2]}$, $\Delta_3,\ldots$ are all involutive. The involutivity of $\Delta_2$ implies that $[\bar{a},\bar{b}_2]\in\Span{\bar{b}_1,\bar{b}_2,[\bar{b}_1,\bar{b}_2]}$ and thus, we actually have $\Delta_2=\Span{\partial_{\bar{u}^1},\bar{b}_1,\bar{b}_2,[\bar{b}_1,\bar{b}_2]}=\Span{\partial_{\bar{u}^1}}+D_1^{(1)}$, which in turn implies that $D_1^{(1)}$ is involutive, \ie $\overline{D}_1=D_1^{(1)}$. Based on these considerations, in the following we explain how to find an input transformation required to generate an input which needs to be prolonged in order to render the system static feedback linearizable. If $n=3$, it can be shown that we can prolong any input $\tilde{u}^1=g(x)^1+m^1_k(x)u^k$ of the system and obtain a static feedback linearizable prolonged system. 
			%			\footnote{Note that we have $D_1^{(1)}+[\bar{a},D_1^{(1)}]=D_1^{(1)}+[a,D_1^{(1)}]$, since $\bar{a}=a+\gamma^1b_1+\gamma^2b_2$ and $D_1^{(1)}$ is involutive. The property $[\bar{a},\bar{b}_1]\notin D_1^{(1)}$ follows from the fact that if $\Dim{\Delta_2}<n+1$, we necessarily have $[a,D_1^{(1)}]\not\subset D_1^{(1)}$ (otherwise, the system contains an autonomous subsystem). Thus, either $[a,\bar{b}_1]\notin D_1^{(1)}$ or $[a,[b_1,b_2]]\notin D_1^{(1)}$. From the Jacobi identity
%			\begin{align*}
%				\begin{aligned}
%					\underbrace{[a,[\bar{b}_1,\bar{b}_2]]}_{\notin D_1^{(1)}}+[\bar{b}_2,[a,\bar{b}_1]]+\overbrace{[\bar{b}_1,\underbrace{[\bar{b}_2,a]}_{\in D_1^{(1)}}]}^{\in D_1^{(1)}}&=0\,,
%				\end{aligned}
%			\end{align*}
%			if follows that $[\bar{b}_2,[a,\bar{b}_1]]\notin D_1^{(1)}$, which because of the involutivity of $D_1^{(1)}$ and $\bar{b}_2\in D_1^{(1)}$ again implies that  $[a,\bar{b}_1]\notin D_1^{(1)}$.}.
			For $n\geq 4$, the direction of the vector field $\bar{b}_2$ is uniquely determined by the condition $[a,\bar{b}_2]\in D_1^{(1)}$ (\ie it can be shown that $[a,\bar{b}_1]\not\in D_1^{(1)}$). The construction is as follows. Calculate $D_1^{(1)}=\Span{b_1,b_2,[b_1,b_2]}$ and find functions $\alpha^1$ and $\alpha^2$ such that $\alpha^1[a,b_1]+\alpha^2[a,b_2]\in D_1^{(1)}$. The functions $\alpha^1$ and $\alpha^2$ are of course only unique up to a multiplicative factor. The vector field $b_c=\alpha^1b_1+\alpha^2b_2$ is then collinear with the vector field $\bar{b}_2$ from above. By applying the input transformation $\tilde{u}^1=\alpha^2u^1-\alpha^1u^2$, $\tilde{u}^2=\beta^1u^1+\beta^2u^2$, with $\beta^1$ and $\beta^2$ chosen such that the transformation is invertible, \ie $det=\alpha^1\beta^1+\alpha^2\beta^2\neq 0$, to the original system, we obtain
			\vspace{-2.5ex}
			\begin{align*}
				\begin{aligned}
					\dot{x}&=%a(x)+\tfrac{1}{det}(b_1(\beta^2\tilde{u}^1+\alpha^1\tilde{u}^2)+b_2(\alpha^2\bar{u}^2-\beta^1\bar{u}^1))\\
				%	&=
					a(x)+\tfrac{1}{det}((\beta^2b_1-\beta^1b_2)\tilde{u}^1+\overbrace{(\alpha^1b_1+\alpha^2b_2)}^{b_c}\tilde{u}^2)\,.
				\end{aligned}
			\end{align*}
			By prolonging the input $\tilde{u}^1$, \ie adding the equation $\dot{\tilde{u}}^1=\tilde{u}^1_1$, we obtain
			\vspace{-1ex}
			\begin{align}\label{eq:d1prolonged}
				\begin{aligned}
					\Sigma_1:&&\dot{x}&=a(x)+\tfrac{1}{det}((\beta^2b_1-\beta^1b_2)\tilde{u}^1+b_c\tilde{u}^2)\\
					\Sigma_2:&&\dot{\tilde{u}}^1&=\tilde{u}^1_1\,,
				\end{aligned}
			\end{align}
			\vspace{-2ex}
			
			\noindent
			with the state $x_p=(x,\tilde{u}^1)$ and the input $u_p=(\tilde{u}^1_1,\tilde{u}^2)$. Calculating the distributions involved in the test for static feedback linearizability of \eqref{eq:d1prolonged}, we obtain the distributions $\Span{\partial_{\tilde{u}^1},b_c}$, $\Span{\partial_{\tilde{u}^1},b_1,b_2,[b_1,b_2]}\,,\ldots$ which all can be shown to be involutive based on the involutivity of the distributions $\Delta_i$ from above, and thus, \eqref{eq:d1prolonged} is static feedback linearizable. The linearizing outputs of \eqref{eq:d1prolonged} are flat outputs with a difference of $d=1$ of the original system $\Sigma_1$. Indeed, let $y$ be a linearizing output of \eqref{eq:d1prolonged}. For the parameterization of the states and inputs of \eqref{eq:d1prolonged} with respect to $y$, we have the diffeomorphism
			\vspace{-1ex}
			\begin{align*}
				\begin{aligned}
					(x,\tilde{u}^1)&=F_{x_p}(y_{[R-1]})\\
					(\tilde{u}^1_1,\tilde{u}^2)&=F_{u_p}(y_{[R]})
				\end{aligned}
			\end{align*}
			\vspace{-2.2ex}
			
			\noindent
			with $\#R=n+1$. It contains the map $(x,\tilde{u}^1)=F_{x_p}(y_{[R-1]})$, $\tilde{u}^2=F_{u_p}^1(y_{[R]})$, which is the parameterization of the state and the input $\tilde{u}$ of the original system with respect to the flat output $y$. It is a submersion with a domain of dimension $\#R+2=n+3$ and a codomain of dimension $n+2$, \ie $d=n+3-(n+2)=1$. (It is immediate that the parameterization of $\tilde{u}^2$ involves the highest derivatives $y_R=(y^1_{r_1},y^2_{r_2})$, otherwise, the original system would be static feedback linearizable, \ie it would have a difference of $d=0$, which contradicts with our assumption $d=1$.) % eventuell ergänzen dass $y=\varphi(x)$ gilt, denn würde $y$ von den Eingängen $\tilde{u}^1$, $\tilde{u}^2$ abhängen, so würde das original system exakt linearisierbar sein, was im Widerspruch zur Annahme $d=1$ steht. 
						
			For $d=2$, it follows from Theorem \ref{thm:d2lin} that there exists an affine input transformation $\bar{u}^j=g(x)^j+m^j_k(x)u^k$ such that the prolonged system
			\vspace{-1ex}
			\begin{align}\label{eq:sfld2}
				&\begin{aligned}
					\dot{x}&=\bar{a}(x)+\bar{b}_1(x)\bar{u}^1+\bar{b}_2(x)\bar{u}^2\\
					\dot{\bar{u}}^1&=\bar{u}^1_1\\
					\dot{\bar{u}}^1_1&=\bar{u}^1_2
				\end{aligned}
			\end{align}
			\vspace{-2ex}
			
			\noindent
			withe the state $(x,\bar{u}^1,\bar{u}^1_1)$ and the input $(\bar{u}^1_2,\bar{u}^2)$ is static feedback linearizable. Thus, the distributions $\Delta_1=\Span{\partial_{\bar{u}^1_1},\bar{b}_2}$, $\Delta_2=\Span{\partial_{\bar{u}^1_1},\partial_{\bar{u}^1},\bar{b}_2,[\bar{a},\bar{b}_2]+\bar{u}^1[\bar{b}_1,\bar{b}_2]}$, $\Delta_3=\Span{\partial_{\bar{u}^1_1},\partial_{\bar{u}^1},\bar{b}_1,\bar{b}_2,[\bar{a},\bar{b}_2]+\bar{u}^1[\bar{b}_1,\bar{b}_2],[\bar{a},[\bar{a},\bar{b}_2]]+\bar{u}^1[\bar{a},[\bar{b}_1,\bar{b}_2]]+\bar{u}^1[\bar{b}_1,[\bar{a},\bar{b}_2]]+(\bar{u}^1)^2[\bar{b}_1,[\bar{b}_1,\bar{b}_2]]+\bar{u}^2[\bar{b}_2,[\bar{a},\bar{b}_2]]+\bar{u}^1\bar{u}^2[\bar{b}_2,[\bar{b}_1,\bar{b}_2]]}$, $\Delta_4,\ldots$ are all involutive. The involutivity of $\Delta_2$ implies that $[\bar{a},\bar{b}_2]\in\Span{\bar{b}_2,[\bar{b}_1,\bar{b}_2]}$, \ie we actually have $\Delta_2=\Span{\partial_{\bar{u}^1_1},\partial_{\bar{u}^1},\bar{b}_2,[\bar{b}_1,\bar{b}_2]}$, which in turn implies that $\bar{b}_2\in\C{D_1^{(1)}}$ (we have $D_1=\Span{\bar{b}_1,\bar{b}_2}$ and $D_1^{(1)}=\Span{\bar{b}_1,\bar{b}_2,[\bar{b}_1,\bar{b}_2]}$). The distribution $\Delta_3$ simplifies to $\Delta_3=\Span{\partial_{\bar{u}^1_1},\partial_{\bar{u}^1},\bar{b}_1,\bar{b}_2,[\bar{b}_1,\bar{b}_2],[\bar{a},[\bar{b}_1,\bar{b}_2]],[\bar{b}_1,[\bar{b}_1,\bar{b}_2]]}$, where the vector fields $[\bar{a},[\bar{b}_1,\bar{b}_2]]$, $[\bar{b}_1,[\bar{b}_1,\bar{b}_2]]$ are linearly dependent$\Mod D_1^{(1)}=\Span{\bar{b}_1,\bar{b}_2,[\bar{b}_1,\bar{b}_2]}$. We have to distinguish between two subcases, namely between $D_1^{(1)}$ being involutive or not. If $D_1^{(1)}$ is non-involutive, we necessarily have $[\bar{b}_1,[\bar{b}_1,\bar{b}_2]]\notin D_1^{(1)}$ and thus $\Dim{D_1^{(2)}}=4$ and $\Dim{\C{D_1^{(1)}}}=1$. Furthermore, we then have $\Delta_3=\Span{\partial_{\bar{u}^1_1},\partial_{\bar{u}^1}}+D_1^{(2)}$, implying that $\overline{D}_1=D_1^{(2)}$. So in this case, the direction of $\bar{b}_2$ is uniquely determined by the condition $\bar{b}_2\in\C{D_1^{(1)}}$. The construction is as follows. Calculate the Cauchy characteristic distribution $\C{D_1^{(1)}}$ and find functions $\alpha^1$ and $\alpha^2$ such that $b_c=\alpha^1b_1+\alpha^2b_2\in\C{D_1^{(1)}}$. The functions $\alpha^1$ and $\alpha^2$ are again only unique up to a multiplicative factor. The vector field $b_c=\alpha^1b_1+\alpha^2b_2$ is then collinear with the vector field $\bar{b}_2$ from above. If $D_1^{(1)}$ is involutive, \ie $\overline{D}_1=D_1^{(1)}$, we obviously have $[\bar{b}_1,[\bar{b}_1,\bar{b}_2]]\in D_1^{(1)}=\Span{\bar{b}_1,\bar{b}_2,[\bar{b}_1,\bar{b}_2]}$ and thus $\Delta_3=\Span{\partial_{\bar{u}^1_1},\partial_{\bar{u}^1},\bar{b}_1,\bar{b}_2,[\bar{b}_1,\bar{b}_2],[\bar{a},[\bar{b}_1,\bar{b}_2]]}$. In this case, the condition $\bar{b}_2\in\C{D_1^{(1)}}$ is not useful since $\C{D_1^{(1)}}=D_1^{(1)}$, but then, analogous to the previously considered case $d=1$, the condition $[\bar{a},\bar{b}_2]\in D_1^{(1)}$  (which follows from $[\bar{a},\bar{b}_2]\in\Span{\bar{b}_2,[\bar{b}_1,\bar{b}_2]}\subset D_1^{(1)}$) uniquely determines the direction of $\bar{b}_2$. The construction of the vector field $b_c=\alpha^1b_1+\alpha^2b_2$ which is collinear with $\bar{b}_2$ is then analogous to the previously considered case $d=1$.
			\begin{remark}\label{rem:functions}
				In the following, we explain how to derive the input transformation required to generate the input which needs to be prolonged twice in order to render the system static feedback linearizable. In contrast to the previously considered case $d=1$, in the case $d=2$, additional functions $\gamma^1,\gamma^2$ and $\delta$ of the state $x$ of the system are involved in the input transformation. We do not discuss their construction. For the proof, only the existence of these functions is of importance and for applying the procedure, we do not have to construct these functions either.
			\end{remark}
		
			In any of the two subcases, \ie independently of $D_1^{(1)}$ being involutive or not, by applying the input transformation $\tilde{u}^1=(\alpha^2u^1-\alpha^1u^2)\delta+\gamma^1$, $\tilde{u}^2=\beta^1u^1+\beta^2u^2+\gamma^2$, with $\beta^1$ and $\beta^2$ chosen such that the transformation is invertible, \ie $det=\alpha^1\beta^1+\alpha^2\beta^2\neq 0$, and functions $\gamma^1,\gamma^2$ and $\delta$ of the state $x$ of the system, to the original system, we obtain
			\vspace{-0.9ex}
			\begin{align*}
				\begin{aligned}
					\dot{x}&=\underbrace{a(x)+\tfrac{\gamma^1}{\delta det}(\beta^2b_1-\beta^1b_2)-\tfrac{\gamma^2}{det}\overbrace{(\alpha^1b_1+\alpha^2b_2)}^{b_c}}_{\tilde{a}}+\\
					&\hspace{3em}\underbrace{\tfrac{1}{\delta det}(\beta^2b_1-\beta^1b_2)}_{\tilde{b}_1}\tilde{u}^1+\tfrac{1}{det}\underbrace{(\alpha^1b_1+\alpha^2b_2)}_{b_c}\tilde{u}^2\,.
				\end{aligned}
			\end{align*}
			\vspace{-1.5ex}
			
			\noindent
			Assume that the functions $\gamma^1$, $\gamma^2$ and $\delta$ are chosen such that $\tilde{a}$ coincides$\Mod\Span{b_c}$ with $\bar{a}$ from above and $\tilde{b}_1$ coincides$\Mod\Span{b_c}$ with $\bar{b}_1$ from above, which is indeed always possible (in fact, we can always choose $\gamma^2=0$). Then, based on the properties of the vector fields $\bar{a}$, $\bar{b}_1$, $\bar{b}_2$ and the distributions $\Delta_i$ form above, it can be shown that we have $[\tilde{a},b_c]\in\Span{b_c,[\tilde{b}_1,b_c]}$ with $\Span{b_c,[\tilde{b}_1,b_c]}$ being involutive and that $[\tilde{a},[\tilde{b}_1,b_c]]$ and $[\tilde{b}_1,[\tilde{b}_1,b_c]]$ are collinear$\Mod D_1^{(1)}$. Prolonging the input $\tilde{u}^1$ two-fold, \ie adding the equations $\dot{\tilde{u}}^1=\bar{u}^1_1$ and $\dot{\tilde{u}}^1_1=\tilde{u}^1_2$ we obtain
			\vspace{-1.3ex}
			\begin{align}\label{eq:d2prolonged}
				\begin{aligned}
					&\Sigma_1:&&\begin{aligned}
						\dot{x}&=\tilde{a}+\tilde{b}_1\tilde{u}^1+\tfrac{1}{det}b_c\tilde{u}^2
					\end{aligned}\\[1ex]
					&\Sigma_2:&&\begin{aligned}
						\dot{\tilde{u}}^1&=\tilde{u}^1_1\\
						\dot{\tilde{u}}^1_1&=\tilde{u}^1_2\,,
					\end{aligned}
				\end{aligned}
			\end{align}
			\vspace{-2.3ex}
			
			\noindent
			with the state $x_p=(x,\tilde{u}^1,\tilde{u}^1_1)$ and the input $u_p=(\tilde{u}^1_2,\tilde{u}^2)$. Calculating the distributions involved in the test for static feedback linearizability of \eqref{eq:d2prolonged}, we then obtain $\Span{\partial_{\tilde{u}^1_1},b_c}$, $\Span{\partial_{\tilde{u}^1_1},\partial_{\tilde{u}^1},b_c,[\tilde{b}_1,b_c]}$, $\Span{\partial_{\tilde{u}^1_1},\partial_{\tilde{u}^1},\tilde{b}_1,b_c,[\tilde{b}_1,b_c],[\tilde{a},[\tilde{b}_1,b_c]],[\tilde{b}_1,[\tilde{b}_1,b_c]]}\,,\ldots$ which all can be shown to be involutive based on the involutivity of the distributions $\Delta_i$ from above. The linearizing outputs of \eqref{eq:d2prolonged} are flat outputs with a difference of $d=2$ of the original system $\Sigma_1$. Indeed, let $y$ be a linearizing output of \eqref{eq:d2prolonged}. For the parameterization of the states and inputs of \eqref{eq:d2prolonged} with respect to $y$, we have the diffeomorphism
			\vspace{-0.6ex}
			\begin{align*}
				\begin{aligned}
					(x,\tilde{u}^1,\tilde{u}^1_1)&=F_{x_p}(y_{[R-1]})\\
					(\tilde{u}^1_2,\tilde{u}^2)&=F_{u_p}(y_{[R]})
				\end{aligned}
			\end{align*}
			\vspace{-1.5ex}
			
			\noindent
			with $\#R=n+2$ (in fact the parameterization of $\tilde{u}^1$ can only involve derivatives of $y$ up to the order $R-2$ due to the structure of \eqref{eq:d2prolonged}). It contains the map
			\vspace{-0.3ex}
			\begin{align*}
				\begin{aligned}
					(x,\tilde{u}^1)&=F_{x_p}^{1,\ldots,n+1}(y_{[R-1]})\\
					\tilde{u}^2&=F_{u_p}^1(y_{[R]})\,,
				\end{aligned}
			\end{align*}
			\vspace{-1.8ex}
			
			\noindent
			which is the parameterization of the state and the input $\tilde{u}$ of the original system with respect to the flat output $y$. It is a submersion with a domain of dimension $\#R+2=n+4$ and a codomain of dimension $n+2$, \ie $d=n+4-(n+2)=2$. (It is immediate that parameterization of $\tilde{u}^2$ involves the highest derivatives $y_R=(y^1_{r_1},y^2_{r_2})$, otherwise, the system would have a difference of $d\leq 1$, which contradicts with our assumption $d=2$.) As already mentioned in Remark \ref{rem:functions}, we do not need to find the appropriate functions $\gamma^1$ and $\delta$. Since $\tilde{u}^1$ with the correctly chosen functions $\gamma^1$ and $\delta$ only depends on $y_{[R-1]}$ (actually $y_{[R-2]}$ due to the structure of \eqref{eq:d2prolonged}) and the functions $\gamma^1$ and $\delta$ only depend on the state $x$ of the original system, which on its own only depends on $y_{[R-1]}$, an arbitrary choice for these functions, \eg $\gamma^1=0$ and $\delta=1$ which would yield the input transformation $\hat{u}^1=\alpha^2u^1-\alpha^1u^2$, $\hat{u}^2=\beta^1u^1+\beta^2u^2$, still yields an input $\hat{u}^1$ which only depends on $y_{[R-1]}$. The flat output $y$ with $d=2$ of the original system is thus a flat output with $d=1$ for the prolonged system
			\vspace{-0.8ex}
			\begin{align}\label{eq:d2onefolProlonged}
				&\begin{aligned}
					\dot{x}&=a(x)+\hat{b}_1(x)\hat{u}^1+\hat{b}_2(x)\hat{u}^2\,,
				\end{aligned}&\begin{aligned}
					\dot{\hat{u}}^1&=\hat{u}^1_1\,,
				\end{aligned}
			\end{align}
			\vspace{-3ex}
			
			\noindent
			\ie \eqref{eq:d2onefolProlonged} has a difference of $d=1$. Conversely, by an analogues reasoning as before based on the parameterizing map of \eqref{eq:d2onefolProlonged} with respect to any flat output with $d=1$ of \eqref{eq:d2onefolProlonged}, it can be shown that any flat output with $d=1$ of \eqref{eq:d2onefolProlonged} is a flat output with $d=2$ of the original system.		
		\hfill\QED
		\vspace{-0.5ex}
		\subsection*{Proof of Lemma \ref{lem:case3}}\noindent
			According to Theorem \ref{thm:d2lin}, the system becomes static feedback linearizable after applying a suitable input transformation $\bar{u}=\Phi_u(x,u)$ and subsequently $d$-fold prolonging the new input $\bar{u}^1$, \ie the prolonged system
			\vspace{-0.8ex}
			\begin{align}\label{eq:d1d2prolonged}
				\begin{aligned}
					\dot{x}&=\bar{f}(x,\bar{u}^1,\bar{u}^2)\\
					\dot{\bar{u}}^1&=\bar{u}^1_1
				\end{aligned}&&\text{or}&&
				\begin{aligned}
					\dot{x}&=\bar{f}(x,\bar{u}^1,\bar{u}^2)\\
					\dot{\bar{u}}^1&=\bar{u}^1_1\\
					\dot{\bar{u}}^1_1&=\bar{u}^1_2\,,
				\end{aligned}
			\end{align}
			\vspace{-2.2ex}
			
			\noindent
			depending on the actual value of $d$, is static feedback linearizable. Analogous to the proof of Lemma \ref{lem:case2}, it can be shown that the linearizing outputs $y$ of \eqref{eq:d1d2prolonged} are flat outputs with a difference of $d=1$ or $d=2$ of the original system and that the parameterization of the input  $\bar{u}^1$ by $y$ involves derivatives up to order $R-d$ only. Independent of the actual value of $d$, the involutivity of $\Span{\partial_{\bar{u}^1_d},\partial_{\bar{u}^2},\partial_{\bar{u}^1_{d-1}},\partial_{\bar{u}^2}\bar{f}^i\partial_{x^i}}$ implies that $\dot{x}=\bar{f}(x,\bar{u}^1,\bar{u}^2)$ is actually of the form $\dot{x}=\bar{a}(x,\bar{u}^1)+b(x,\bar{u}^1)h(x,\bar{u}^1,\bar{u}^2)$ and introducing $\tilde{u}^2=h(x,\bar{u}^1,\bar{u}^2)$ thus results in the PAI representation
			\vspace{-0.4ex}
			\begin{align}\label{eq:specialPAI}
				\begin{aligned}
					\dot{x}&=a(x,\bar{u}^1)+b(x,\bar{u}^1)\tilde{u}^2\,.
				\end{aligned}
			\end{align}
			\vspace{-3.2ex}
			
			\noindent
			This PAI representation is special, since the parameterization of the non-affine occurring input $\bar{u}^1$ involves derivatives of $y$ up to order $R-d$ only. For every PAI representation which is equivalent to \eqref{eq:specialPAI} via a PAI form preserving transformation \eqref{eq:paiPreservingTransformation}, we still have that the non-affine input only depends on $y_{[R-1]}$, since in such transformations the non-affine input is only combined with the state $x$ of the system and we have $x=F_x(y_{[R-1]})$. Recall that there exist at most two fundamentally different PAI representations of a system, all others are equivalent to one of those two by such a PAI form preserving transformation \eqref{eq:paiPreservingTransformation}. Therefore, given any two non-equivalent PAI representations of the system (which can be derived systematically, provided they exist), one of them is equivalent to the special PAI representation \eqref{eq:specialPAI} via a transformation of the form \eqref{eq:paiPreservingTransformation}. The non-affine inputs of any two non-equivalent PAI representations are thus the candidates for an input whose parameterization involves derivatives of $y$ up to order $R-1$ only. Similar as in the proof of Lemma \ref{lem:case2}, it can then be shown that at least one of the possibly two AI-systems obtained by prolonging the corresponding non-affine inputs has a difference of $d\leq 1$ and that furthermore the linearizing outputs or flat outputs with $d=1$ of the AI-systems are flat outputs with $d=1$ or $d=2$ of the original system.\hfill\QED
			\vspace{-0.5ex}
		\subsection*{Proof of Theorem \ref{thm:d2lin}}\noindent
			The proof of Theorem \ref{thm:d2lin} is based on the following results.
			\begin{lemma}\label{lem:d2xu}
				A system \eqref{eq:sys} with $d\leq 2$ is $(x,u)$-flat, an AI-system \eqref{eq:sysAI} with $d\leq 2$ is $x$-flat.
			\end{lemma}
		
			Due to space limitations, we do not provide a proof of this lemma here. The second part of the lemma, \ie that AI-systems with $d\leq 2$ are $x$-flat, can also be found in \cite{NicolauRespondek:2016-2}, Proposition 1.
%			As a consequence, such systems can be rendered static feedback linearizable by $d$-fold prolonging a suitably chosen (new) control, as the following theorem asserts.
			\begin{theorem}\label{thm:linearizationByProlongations}
				Every $(x,u)$-flat system with two inputs can be rendered static feedback linearizable by $d=\#R-n$ fold prolonging a suitably chosen (new) input.
				%\footnote{Here, $d$ denotes the difference of a not necessarily minimal $(x,u)$-flat output of the system. To every $(x,u)$-flat output $y$ of the system, there exists a static feedback linearizable prolonged system such that $y$ is a linearizing output of the prolonged system.}.
			\end{theorem}
		
			A proof of this theorem can be found in \cite{GstottnerKolarSchoberl:2020b}. Since according to Lemma \ref{lem:d2xu}, systems with $d\leq 2$ are $(x,u)$-flat, Theorem \ref{thm:linearizationByProlongations} always applies to these systems and in turn, such systems can always be rendered static feedback linearizable by prolongations of a suitably chosen (new) input, which shows the first part of Theorem \ref{thm:d2lin}. From the proof of Theorem \ref{thm:linearizationByProlongations} in \cite{GstottnerKolarSchoberl:2020b}, it follows that the input transformation which has to be applied in order to generate the required input which needs to be prolonged, is given by $\bar{u}^1=\Lie_f^{k_j}\varphi^j(x,u)$, $j=1$ or $j=2$ where $k_j\geq 0$ denotes the relative degree of the component $y^j=\varphi^j(x,u)$ of the flat output, \ie $\Lie_f^{k_j-1}\varphi^j=\Lie_f^{k_j-1}\varphi^j(x)$ and $\Lie_f^{k_j}\varphi^j=\Lie_f^{k_j}\varphi^j(x,u)$. Since according to Lemma \ref{lem:d2xu} an AI-system with $d\leq 2$ is $x$-flat, \ie $y^j=\varphi^j(x)$, there always exists an affine input transformation which generates the required input since $\bar{u}^1=\Lie_{(a^i+u^1b_1^i+u^2b_2^i)\partial_{x^i}}^{k_j}\varphi^j(x)$ is an affine input transformation.\hfill\QED
	\bibliographystyle{IEEEtran} 
	\bibliography{IEEEabrv,Bibliography}

\end{document}